\documentclass{ifacconf}

\usepackage{natbib}        
\usepackage{amsfonts}
\usepackage{amsmath, bbm}
\usepackage{amssymb}
\usepackage{amsxtra}
\usepackage{psfrag}
\usepackage{ifthen}
\usepackage[dvipsnames]{xcolor}
\usepackage{color}
\usepackage{datetime}
\usepackage[utf8]{inputenc}
\usepackage[english]{babel}
\usepackage{multicol}
\usepackage{mathtools}
\usepackage{graphicx}      
\usepackage{epstopdf}
\usepackage{url}

\usepackage{enumitem}

\newboolean{showcomments}
\setboolean{showcomments}{true}

\newcommand{\khaled}[1]{  \ifthenelse{\boolean{showcomments}}
{ \textcolor{blue}{(Khaled says:  #1)}} {}  }
\newcommand{\bose}[1]{\ifthenelse{\boolean{showcomments}}
{ \textcolor{red}{(Bose says:  #1)}}{}}
\newcommand{\tamer}[1]{\ifthenelse{\boolean{showcomments}}
{ \textcolor{red}{(Tamer says:  #1)}}{}}

\newcommand{\addcite}{{\textcolor{brown}{(addcite) }}}

\newcommand{\beq}{\begin{equation}}
\newcommand{\eeq}{\end{equation}}
\newcommand{\beqa}{\begin{eqnarray}}
\newcommand{\eeqa}{\end{eqnarray}}
\newcommand{\beqan}{\begin{eqnarray*}}
\newcommand{\eeqan}{\end{eqnarray*}}

\newcommand{\E}{\mathbb{E} }

\renewcommand{\prob}{\mathbb{P}}
\newcommand{\ind}[1]{\mathbbm{1}_{\lb#1\rb}}

\newcommand{\Gset}{\mathbb{G}}

\renewcommand{\Nset}{\mathbb{N}}

\newcommand{\Pset}{\mathbb{P}}

\renewcommand{\Rset}{\mathbf{R}}

\newcommand{\Acal}{{\cal A}}

\newcommand{\Ecal}{{\cal E}}
\newcommand{\Fcal}{{\cal F}}
\newcommand{\Gcal}{{\cal G}}

\newcommand{\Ncal}{{\cal N}}

\newcommand{\Zcal}{{\cal Z}}

\newcommand{\Gfk}{{\mathfrak{G}}}

\newcommand{\Rfk}{{\mathfrak{R}}}

\newcommand{\bone}{\mathbf{1}}

\renewcommand{\[}{\left[}
\renewcommand{\]}{\right]}
\newcommand{\lb}{\left\{}
\newcommand{\rb}{\right\}}
\renewcommand{\(}{\left(}
\renewcommand{\)}{\right)}

\newcounter{l1}
\newcounter{l2}
\newcounter{l3}
\newcommand{\bdotlist}{\begin{list}{$\bullet$}{}}
\newcommand{\bboxlist}{\begin{list}{$\Box$}{}}
\newcommand{\bbboxlist}{\begin{list}{\raisebox{.005in}{{\tiny
$\blacksquare$ \ \ }}}{}}
\newcommand{\bdashlist}{\begin{list}{$-$}{} }
\newcommand{\blist}{\begin{list}{}{} }
\newcommand{\barablist}{\begin{list}{\arabic{l1}}{\usecounter{l1}}}
\newcommand{\balphlist}{\begin{list}{(\alph{l2})}{\usecounter{l2}}}
\newcommand{\bAlphlist}{\begin{list}{\Alph{l2}.}{\usecounter{l2}}}
\newcommand{\bdiamlist}{\begin{list}{$\diamond$}{}}
\newcommand{\bromalist}{\begin{list}{(\roman{l3})}{\usecounter{l3}}}

\renewcommand{\Nset}{{\mathfrak{N}}}
\renewcommand{\Gset}{{\mathfrak{G}}}

\renewcommand{\Pset}{{\mathfrak{P}}}

\renewcommand{\Rset}{{\mathbf{R}}}
\renewcommand{\prob}{{\mathbf{P}}}

\begin{document}
\begin{frontmatter}

\title{Cash-settled options for wholesale electricity markets} 


\author{Khaled Alshehri,} 
\author{Subhonmesh Bose,}   
\author{and Tamer Ba\c{s}ar}

\address{Department of Electrical and Computer Engineering, 
   University of Illinois at Urbana-Champaign, Urbana, IL 61801 USA  (e-mail: \{kalsheh2,boses,basar1\}@illinois.edu.)}
   
\begin{abstract}:	
Wholesale electricity market designs in practice do not provide the market participants with adequate mechanisms to hedge their financial risks. Demanders and suppliers will likely face even greater risks with the deepening penetration of variable renewable resources like wind and solar. This paper explores the design of a centralized cash-settled call option market to mitigate such risks. A cash-settled call option is a financial instrument that allows its holder the right to claim a monetary reward equal to the positive difference between the real-time price of an underlying commodity and a pre-negotiated strike price for an upfront fee. Through an example, we illustrate that a bilateral call option can reduce the payment volatility of market participants. Then, we design a centralized clearing mechanism for call options that generalizes the bilateral trade. We illustrate through an example how the centralized clearing mechanism generalizes the bilateral trade. Finally, the effect of risk preference of the market participants, as well as some generalizations are discussed. 
\end{abstract}
\begin{keyword}
Electricity Markets, Call Options, Mechanism Design, Stackelberg Equilibrium.
\end{keyword}

\end{frontmatter}




\section{Introduction}
Various states in the U.S. and countries around the world have adopted aggressive targets for the integration of renewable energy resources. Wind and solar energy are two of the most prominent resources. The inherent variability of these resources makes it difficult to maintain the balance of demand and supply of power at all times. By variable, we mean they are uncertain (errors in day-ahead forecasts are significantly higher than those in bulk power demand), intermittent (shows large ramps over short time horizons), and non-dispatchable (output cannot be varied on command). See \cite{Bird} for a comprehensive discussion on the challenges of renewable integration.

Energy is typically procured in advance to meet the demand requirements. Forward planning is necessary since many generators -- such as the ones based on nuclear technology or coal -- cannot alter their outputs arbitrarily fast to track demand requirements; some lead time is necessary. In its simplest abstraction, one can model the system operation to proceed in two stages: a forward stage, conducted a day or a few hours in advance, and the real-time stage. Roughly, the forward stage optimizes the dispatch against a forecast of the demand and supply conditions at real-time. The impending deviation from such forecasts are then balanced in real-time. While demand forecasts even a day in advance are within 1-3\% accuracy, the same forecasts in the availability of variable renewable resources can be significantly higher; they can be as high as 12\%\footnote{Some promising forecasting techniques have been known to reduce the forecast error further to 6-8\% over large geographical regions.} For more details on the statistics, see \cite{Bird}. Variability in supply from resources like wind and solar exposes market participants to increased financial risks. The forward market design in practice does not allow participants in the wholesale market to adequately hedge their financial risks. This paper proposes a financial instrument for the same. The deepening penetration of variable renewable supply will increase the volatility in payments to market participants and hence, increase the financial risks borne by market participants (see \cite{wang, marketevo}). One needs to better design the financial instruments to mitigate such risks. 

Electricity market participants engage in trading financial derivatives, i.e., instruments that derive their values based on the prices in the wholesale  market. Traded financial derivatives in practice include electricity forwards, futures, swaps, and options. See \cite{der, swingsurvey, kluge} and the references therein. Some are traded on an exchange and many are traded bilaterally.

In this paper, we consider how an intermediary (called the `market maker') can convene a financial market for trading in {\it{cash-settled call options}}, and how such a financial market can reduce payment volatilities of wholesale electricity market participants. Upon buying one unit of a cash-settled call option at a negotiated option price, the buyer is entitled to receive a cash payment equal to the real-time price of a commodity (that is electricity in our case) less the negotiated strike price. In order to apply the market design to electricity markets, we adopt an economic dispatch and pricing model in Section \ref{sec:model} to provide us with a real-time (spot) price of electricity. Then, we use that model on a stylized example in Section \ref{sec:motivate} and illustrate that a bilateral trade in cash settled call-options between a renewable power producer and a dispatchable peaker power plant can lower their respective payment volatilities. We recognize that engaging in multiple bilateral option trades on a daily basis will likely be difficult in a practical electricity market setting, and will adversely affect the liquidity of such trades.\footnote{As such, markets for financial derivatives associated with electricity markets have been known to suffer from low liquidity. For example, see \cite{liquidity}.} The remedy we offer is a centralized market clearing mechanism for call options in Section \ref{sec:centralized} for its use among electricity market participants. Such a mechanism will make option trading more viable in practice and attractive to market participants. We delineate the salient features and discuss possible generalizations of our design in Section \ref{sec:risk}, and and we conclude in Section \ref{sec:conc}.

Our proposed mechanism is compatible with alternate wholesale electricity market designs, as given in  \cite{Wong, singleset, Bouffard1, Bouffard2, bose}. Even with different designs, market participants can have incentives to strategize their actions in the electricity market and the option trade together.\footnote{We refer to \cite{ledgerwood, chiaraHogan} for discussions on how market participants can strategize their actions across the electricity markets and their associated financial instruments.}Such interactions can adversely affect the market outcomes. However, we relegate such considerations for future work.

\subsubsection{Notation:}
We let $\Rset$ denote the set of real numbers, and $\Rset_+$ (resp. $\Rset_{++}$) denote the set of nonnegative (resp. positive) numbers. For $z\in\Rset$, we let $z^+ := \max\{z, 0\}$.
We let $\E [Z]$ denote the expectation of a random variable $Z$. For any set $\Zcal$, we denote its cardinality by $|\Zcal|$. 
For an event $\Ecal$, we denote its probability by $\prob\{\Ecal\}$ for a suitably defined probability measure $\prob$. The indicator function for an event $\Ecal$ is given by 
\begin{align*} \ind{\Ecal} := \begin{cases} 1, & \text{if } \Ecal \text{ occurs},\\ 0, & \text{otherwise}.\end{cases}\end{align*}
In any optimization problem, a decision variable $x$ at optimality is denoted by $x^*$.


\section{Describing the marketplace}  \label{sec:model}
\vspace{-0.1in}

The wholesale electricity market is comprised of consumers and producers of electricity. The consumers in this market are the load-serving entities that represent the retail customers they serve within their geographical footprint. Examples of such load-serving entities are the utility companies and retail aggregators. Bulk power generators are the producers in this market. We distinguish between two sets of generators. The first type is a \emph{dispatchable generator} that can alter its output within its capabilities on command. Such generators are fuel based; e.g., they run on nuclear technology, or fossil fuels like coal or natural gas, or dispatchable renewable resources such as biomass or hydro power. The second type is a \emph{variable renewable power producer}. Its available capacity of production depends on an intermittent resource like wind or solar irradiance. The system operator, denoted by $SO$, implements a centralized market mechanism that determines the production and consumption of each market participant and their compensations. It does so in a way that balances demand with supply, and the power injections across the grid induce feasible power flows over the transmission lines. Most electricity markets in the United States have a locational marginal pricing based compensation scheme. In this paper, we ignore the transmission constraints of the grid and hence, describe an electricity market with a marginal pricing scheme. Generalizing this work to the case with a network is left for future work.

\subsubsection{Modeling uncertainty in supply:}
To model uncertainty in supply conditions, we consider a two-period market model as follows. Let $t=0$ denote the ex-ante stage, prior to the uncertainty being realized. At this stage, one only has forecasts of the uncertain parameters. The uncertainty is realized at $t=1$, the ex-post stage. One can identify $t=0$ as the \emph{day-ahead stage} and $t=1$ as the \emph{real-time stage} in electricity market operations. Let $(\Omega, \Fcal, \prob)$ denote the probability space, describing the uncertainty. Here, $\Omega$ is the collection of possible scenarios at $t=1$ (which could be uncountably infinite) \footnote{In this paper, we use ``distribution" and ``measure" interchangeably, as appropriate.}, $\Fcal$ is a suitable $\sigma$-algebra over $\Omega$, and $\prob$ is a probability distribution over $\Omega$.
We assume that all market participants know $\prob$.

\subsubsection{Modeling the market participants:}
Let $d$ denote the aggregate inflexible demand that is accurately known a day in advance\footnote{Day-ahead demand forecasts in practice are typically  quite accurate. Notwithstanding the availability of such forecasts, our work can be extended to account for demand uncertainties.}.
Let $\Gfk$ and $\Rfk$ denote the collection of dispatchable generators and variable renewable power producers, respectively. We model their individual capabilities as follows. 

\begin{itemize}[leftmargin=*]
\item Let each dispatchable generator $g \in \Gfk$ produce $x^\omega_g$ in scenario $\omega\in\Omega$. We model its ramping capability by letting $ \lvert x_g^\omega  - x_g^0 \rvert \leq \ell_g,$
where $x_g^0$ is a generator set point that is decided ex-ante, and $\ell_g$ is the ramping limit. Let the installed capacity of generator $g$ be ${x}^{\text{cap}}_g$, and hence $x_g^\omega \in [0, {x}^{\text{cap}}_g]$. Its cost of production is given by the smooth convex increasing map $c_g : [0, {x}^{\text{cap}}_g] \to \Rset_+$.

\item Each variable renewable power producer $r \in \Rfk$ produces $x_r^\omega$ in scenario $\omega\in\Omega$. It has no ramping limitations, but its available production capacity is random, and we have $ x_r^\omega \in [0, \overline{x}_r^\omega] \subseteq [0, {x}^{\text{cap}}_r].$
That is, $\overline{x}_r^\omega$ denotes the random available capacity of production, and $ {x}^{\text{cap}}_r$ denotes the installed capacity for $r$. Similar to a dispatchable generator, the cost of production for $r$ is given by the smooth convex increasing map $c_r : [0, {x}^{\text{cap}}_r] \to \Rset_+$.
\end{itemize}
We call a vector comprised of $x_g$ for each $g \in \Gfk$ and $x_r$ for each $r\in\Rfk$ a \emph{dispatch}.

\subsubsection{Conventional dispatch and pricing model:}

The $SO$ balances demand and supply of power in each scenario. It determines the dispatch and the compensations of all market participants. In the remainder of this section, we describe the so-called conventional dispatch and pricing scheme. This market mechanism serves as a useful benchmark for electricity market designs under uncertainty, e.g., in \cite{Morales1, Morales2}.

We assume that the SO knows $c_g, x_g^{\text{cap}}$ for each $g\in\Gfk$ and $x_r, x_r^{\text{cap}}, \overline{x}_r$ for each $r \in \Rfk$. In practice, the cost functions are derived from supply offers from the generators. The market participants, in general, may have incentives to misrepresent their cost functions. Analyzing the effects of such strategic behavior is beyond the scope of this paper.

%

\subsubsection{The day-ahead stage:}

The SO computes a forward dispatch against a point forecast of all uncertain parameters. In particular, the SO replaces the random available capacity $\overline{x}^\omega_r$ by a certainty surrogate $\overline{x}^{\text{CE}}_r \in [0, \overline{x}^{\text{cap}}_r] $ for each $r \in \Rfk$, and computes the forward dispatch by solving
\begin{equation*}
\begin{alignedat}{5}
&{\text{minimize}}   \ \ \ 
	& & \sum_{g \in \Gfk} c_g(X_g) + \sum_{r \in \Rfk} c_r(X_r), \\
& \text{subject to}  & & \sum_{g \in \Gfk} X_g + \sum_{r \in \Rfk} X_r = d, \\
&&& X_g \in [0, {x}_g^{\text{cap}}], \ \ X_r \in [0, \overline{x}^{\text{CE}}_r],\\
&&& \text{for each } g \in \Gfk, \ r \in \Rfk,
\end{alignedat}
\label{eq:DAM}
\end{equation*}
over $X_g \in \Rset, g\in \Gfk$, and $X_r \in \Rset, r \in \Rfk$. A popular surrogate\footnote{See \cite{Morales2} for an alternate certainty surrogate.} is given by $ \overline{x}^{\text{CE}}_r := \E[\overline{x}_r^\omega]$.

The forward price is given by the optimal Lagrange multiplier of the energy balance constraint. Denoting this price by $P^*$, generator $g \in \Gfk$ is paid $P^* X_G^*$, while producer $r \in \Rfk$ is paid $P^* X_r^*$. Aggregate consumer pays $P^* d$.

%

%
%

\subsubsection{At real-time:}
Scenario $\omega$ is realized, and the SO solves
\begin{equation*}
\begin{alignedat}{5}
&{\text{minimize}}   \ \ \ 
	& & \sum_{g \in \Gfk} c_g(x_g^\omega) + \sum_{r \in \Rfk} c_r(x_r^\omega), \\
& \text{subject to}  & & \sum_{g \in \Gfk} x_g^\omega + \sum_{r \in \Rfk} x_r^\omega = d, \\
&&& x_g^\omega \in [0, {x}_g^{\text{cap}}], \lvert x_g^\omega  - X_g^* \rvert \leq \ell_g,\\
&&& x_r^\omega \in [0, \overline{x}^\omega],\ \ \text{for each } g \in \Gfk, \ r \in \Rfk,
\end{alignedat}
\label{eq:RTM}
\end{equation*}
over $x_g^\omega \in \Rset, g\in \Gfk$ and $x_r^\omega \in \Rset, r \in \Rset$. The real-time (or spot) price is again defined by the optimal Lagrange multiplier of the energy balance constraint, and is denoted by $p^{\omega,*} \in\Rset_+$. Note that the optimal $X_g^*$ computed at $t=0$ defines the generator set-points $x_g^0$ for each generator $g\in\Gfk$. Generator $g \in \Gfk$ is paid $p^{\omega,*} \left( x_g^{\omega,*}- X_g^* \right)$, while producer $r \in \Rfk$ is paid $p^{\omega,*} \left( x_r^{\omega,*}- X_r^* \right)$. The aggregate consumer does not have any real-time payments, since there is no deviation in the demand.

 
The total payments to each participant is the sum of her forward and the real-time payments. Call these payments as $\pi_g^\omega$ for each $g \in \Gfk$ and $\pi_r^\omega$ for each $r \in \Rfk$ in scenario $\omega$. 
\vspace{-0.1in}

We remark that the conventional dispatch model generally defines a suboptimal forward dispatch in the sense that the generator set-points are  \emph{not} optimized to minimize the expected aggregate costs of production. Several authors have advocated a so-called stochastic economic dispatch model, wherein the forward set-points are optimized against the expected real-time cost of balancing. See \cite{Wong, singleset, Bouffard1, Bouffard2, bose}. Our financial market is designed to work in parallel to any electricity market. Different designs of the latter can be accommodated. We adopt the conventional model to make our discussion and examples concrete.



\section{An example with bilateral cash-settled call option} \label{sec:motivate}
\vspace{-0.1in}

In this section, we study a simple power system example and illustrate that a bilateral cash-settled call option can reduce the variance of the payments and even mitigate the risks of negative payments or financial losses to market participants. This example will motivate our study of a centralized call option market in the next section.

Consider an example power system (adopted from \cite{bose}) with two dispatchable generators and a single variable renewable power producer: $ \Gset := \{ B, P \}, \quad \text{and} \quad \mathfrak{R} : = \{ W \}$,
where $B$ is a base-load generator, $P$ is a peaker power plant, and $W$ is a wind power producer. 

Let $x_B^{\text{cap}} = x_P^{\text{cap}} = \infty, \quad \text{and} \quad \ell_B = 0, \ \ell_P= \infty.$
Therefore, $B$ and $P$ have infinite generation capacities. $B$ cannot alter its output in real-time from its forward set-point, and $P$ has no ramping limitations. Suppose $B$ and $P$ both have linear costs of production. $B$ has a unit marginal cost, while $P$ has a marginal cost of $1/{\rho}$. Let $\rho \in (0,1]$, i.e., $P$ is more expensive than $B$.

To model the uncertainty in available wind, we let
$$ \Omega := [\mu - \sqrt{3}\sigma, \mu + \sqrt{3}\sigma] \subset \Rset_+,$$
and take $\prob$ to be the uniform distribution over $\Omega$. One can verify that $\mu$ and $\sigma^2$ define the mean and the variance of the distribution, respectively. Then, the available wind capacity in scenario $\omega$ is encoded as $\overline{x}_r^\omega = \omega$. Further, assume that $W$ produces power at zero cost, and that $d \geq \mu + \sqrt{3} \sigma$.

This stylized example is a caricature of electricity markets with deepening penetration of variable renewable supply. Base-load generators, specifically those based on nuclear technology, have limited ramping capabilities. Peaker power plants based on internal combustion engines can quickly ramp up their power outputs. Utilizing them to balance variability, however, is costly. Finally, demand is largely inflexible and can be predicted with high accuracy. 

In what follows, we analyze the effect of a bilateral call option on this market example. Insights from this example will prove useful in the balance of the paper.

\subsubsection{Conventional dispatch and pricing for the example power system:}
\vspace{-0.05in}

The conventional dispatch model yields the following forward dispatch and forward price.
$$ X_B^* = d - \mu, \ X_P^* = 0, \ X_W^* = \mu, \ P^* = 1. $$
In scenario $\omega$, the real-time dispatch is given by
$$x_B^{\omega,*} = d - \mu, \ x_W^{\omega,*} = \min\{\omega, \mu\}, \ x_P^{\omega,*} = (\mu - \omega)^+,$$
and the real-time price is given by
\begin{equation}
p^{\omega,*} = (1/\rho) \ind{\omega \in \[ \mu-\sqrt{3}\sigma, \mu \]},
\label{spot} 
\end{equation}
where $\ind{\cdot }$ denotes the indicator function.
The dispatch and the prices yield the following payments in scenario $\omega$:
\begin{align*}
\pi_B^\omega = d-\mu, 
\quad \pi_P^\omega = (\mu-\omega)^+/\rho, \quad
\pi_W^\omega = \mu - ( \mu - \omega)^+/\rho. 
\end{align*}

Next, we consider a bilateral cash-settled option between $P$ and $W$. We will show that  such an option can reduce the volatility of the payments to $W$ and $P$. Also of interest is $W$'s payment, when $\rho < \sqrt{3} \frac{\sigma}{\mu}$ and the realized available wind $\omega$ is in $\Omega_0^-$, where
\begin{align}
\Omega_0^- := \lb \omega \ : \ \mu - \sqrt{3} \sigma \leq \omega < \mu (1 - \rho) \rb.
\label{eq:Omega0}
\end{align}
For such a scenario $\omega$, we have $\pi_W^\omega < 0$. That is, $W$ incurs a loss for each $\omega \in \Omega_0^-$. As we will illustrate, $W$ can avoid incurring such a loss in a subset of the scenarios by trading in an option with $P$. 


\subsection{Bilateral cash-settled call option between $W$ and $P$}
\vspace{-0.05in}

Suppose $W$ buys cash-settled call options from $P$. One unit of such an  option entitles $W$ to a claim of a cash payment of the difference between the real-time price and a pre-negotiated strike price from $P$. We will discuss how such an option  reduces their volatility of payments.

We model the bilateral option trade between $P$ and $W$ as a robust Stackelberg game (see \cite{basar}) $\Gcal$ as follows. Right after the day ahead market is settled at $t=0$, $P$ announces an option price $q \in \Rset_+$ and a strike price $K \in \Rset_+$ for the call option it sells. Then, $W$ responds by purchasing $\Delta \in \Rset_+$ options.\footnote{For simplicity, we allow the possibility of a fractional number of options being bought.} This option entitles $W$ to a cash payment of $\left( p^{\omega,*} - K \right)^+\Delta$ from $P$ in scenario $\omega$. The option costs $W$ a fee of $q\Delta$. Assume that there is an exogenously defined cap of $\sqrt{3}\sigma$ on the amount of option $W$ can buy from $P$. The cap equals the maximum loss that $W$ can incur from the electricity market in real-time.

Then, the total payments to $P$ and $W$ in scenario $\omega$ are given by
\begin{align}
\Pi^\omega_W(q, K, \Delta) &:= \pi_W^\omega - q \Delta + \left( p^{\omega,*} - K \right)^+\Delta,\label{eq:PiW}\\
\Pi^\omega_P(q, K, \Delta) &:= \pi_P^\omega + q\Delta - \left( p^{\omega,*} - K \right)^+\Delta, \label{eq:PiP}
\end{align}
respectively, when they agree on the triple $\(q,K,\Delta\)$. 
Assume that $W$ is risk-neutral and has the correct conjectures on the real-time prices. Then, the perceived payoff for $W$ in the day-ahead stage is given by $\E \left[\Pi^\omega_W(q, K, \Delta)\right]$. If $W$ believes in a different price distribution, the perceived payoff for $W$ can be modified accordingly. Moreover, if $W$ is risk-averse, one can replace $\E$ with a suitable risk-functional. Similar considerations extend to the perceived payoff for $P$.

The possible outcomes of the option trade are identified as the set of Stackelberg equilibria (SE) of $\Gcal$. Precisely, we say $(q^*, K^*, \Delta^*(q^*, K^*))$ constitutes a Stackelberg equilibrium of $\Gcal$, if 
$$\E\left[\Pi_P^\omega(q^*, K^*, \Delta^*(q^*, K^*))\right] \geq \E\left[\Pi_P^\omega(q, K, \Delta^*(q, K))\right]$$
where $\Delta^*:\Rset^2_+\rightarrow[0,\sqrt{3}\sigma]$ is the best response of $W$ to the prices $(q,K)\in\Rset^2_+$ announced by $P$. For a given $(q,K)$, the best response $\Delta^*$ satisfies $$\E\left[\Pi_W^\omega(q, K, \Delta^*(q, K))\right] \geq \E\left[\Pi_W^\omega(q, K, \Delta(q, K))\right],$$ for all $\Delta:\Rset^2_+\rightarrow[0,\sqrt{3}\sigma]$.
 In the following result, we characterize all Stackelberg equilibria of $\Gcal$. In presenting the result, we make use of the notation ${\sf var}[z^\omega]$ to denote the variance of a real valued map $z:\Omega \to \Rset$.

%
\begin{prop}
\label{prop:bilat}
The Stackelberg equilibria of $\Gcal$ are given by $\Ncal_1 \cup \Ncal_2$, where
\begin{align*}
\Ncal_1 &:= \left\{ (q, K, \Delta) \ \vert \ (q,K) \in \Rset_+^2, \ \Delta:\Rset^2_+ \to [0, \sqrt{3}\sigma], \right.\\
&\qquad\quad \left. 2q + K > 1/\rho , \ \Delta(q', K') = 0 \ \forall \ (q', K') \in\Rset^2_+  \right\},\\ 
\Ncal_2 &:= \left\{ (q, K, \Delta) \ \vert \ (q,K) \in \Rset_+^2, \ \Delta:\Rset^2_+ \to [0, \sqrt{3}\sigma], \right.\\ &\qquad\quad \left. 2q + K = 1/\rho \right\}.
\end{align*}
Moreover, for any $(q^*, K^*, \Delta^*(q^*, K^*)=\sqrt{3} \sigma) \in \Ncal_2$, the variances of the payments satisfy
$$
{\sf var}\left[\Pi^\omega_i(q^*, K^*, \Delta^*(q^*, K^*))\right]
- {\sf var}\left[\pi^\omega_i\right] 
= - 3{{K^*}}{\sigma^2} /2< 0,
$$
for each $i \in \lb W, P \rb$.
\end{prop}
%
%
Proof of this Proposition can be found in Appendix A.

We remark that it can be verified that volatility in payments -measured in terms of the respective variance- will decrease for any nontrivial SE. However, our choice of $\Delta^*(q^*,K^*)=\sqrt{3}\sigma$ reveals the main insights without loss of generality. $\Ncal_1$ describes the degenerate case, where $P$ and $W$ essentially do not participate in the option market. 

Proposition 1 implies that $K^* = 1/\rho - 2q^*$ whenever a trade occurs. Thus, for a given option price, the reduction in variance increases as $\rho$ decreases. Said differently, both participants gain more in terms of reduction in volatility as the peaker plant becomes more costly, leading to the possibility of price spikes in the real-time electricity market. Further, the mentioned reduction in variance increases with $\sigma$. Again, the participants stand to gain more (in terms of volatility in payments) from the option trade as the available wind becomes more uncertain.

Finally, consider the collection of scenarios $\omega$ where $W$ suffers a financial loss in the energy market, i.e., $\pi_W(\omega) < 0$. Recall that this happens when $\rho < \sqrt{3} \frac{\sigma}{\mu}$ and $\omega \in \Omega_0^-$ as defined in \eqref{eq:Omega0}.
When $W$ and $P$ participate in the bilateral option trade, then $W$ suffers a loss only when $\omega \in \Omega^-$, where
\begin{equation}\Omega^-:=\{\omega: \mu - \sqrt{3} \sigma \leq \omega < \mu (1 - \rho)-{\rho} q^*\sqrt{3}\sigma\}. \label{eq:subset} \end{equation}
Note that $\Omega^-$ is a strict subset of $\Omega_0^-$. Thus, $W$ is less exposed to negative payments with the option trade.
\subsubsection{Limitations of bilateral trading:} 
\vspace{-0.1in}

Our stylized example reveals the benefits of a bilateral trade in call options. In a wholesale market with $\Gfk$ and $\Rfk$ describing the set of dispatchable and variable generators, respectively, one can conceive of $\lvert \Gfk \rvert \cdot \lvert \Rfk \rvert $ bilateral trade agreements. It is difficult to convene and settle such a large number of trades on a regular basis. To circumvent this difficulty, we propose a central clearing mechanism for cash-settled call options among the market participants.

\section{Centralized market clearing for cash-settled call options}
 \label{sec:centralized}
 \vspace{-0.1in}

Consider an intermediary who acts as an aggregate buyer of such options for a collection of option sellers, and acts as a seller of such options for a collection of option buyers. Call this intermediary a \emph{market maker}, denoted by ${M}$. An established financial institution or even the ${SO}$ can fulfill the role of such a market maker. Let $\Rfk$ be the set of buyers and $\Gfk$ be the set of sellers of the call options. The option market proceeds as follows.

\subsubsection{Day-ahead stage: }
\vspace{-0.1in}
\begin{itemize}[leftmargin=*]

\item $M$ broadcasts a set of allowable trades $\Acal_0$, given by
$$\Acal_0 := \[0, \overline{q}\] \times \[0, \overline{K}\] \times \[0, \overline{\Delta}\] \subset \Rset^3_{+},$$
to all the market participants $\Gfk \cup \Rfk$. 

\item Each $i \in \Gfk \cup \Rfk$ submits an acceptable set of option trades, denoted by $\Acal_i \subseteq \Acal_0$.\footnote{One can fix a parametric description, where market participants communicate their choices of parameters.}

\item $M$ correctly conjectures the real-time prices $p^{\omega,*}$ in each scenario $\omega$, and solves the following stochastic optimization problem to clear the option market.%
\begin{align}
\label{eq:DA.opt}
\begin{alignedat}{5}
&{\text{maximize}}   \ \ \ \E[f^\omega], \\
& \text{subject to} \\
& \qquad  \sum_{g \in \Gfk} \Delta_g = \sum_{r \in \Rfk} \Delta_r, \\
& \qquad (q_g, K_g, \Delta_g)\in\Acal_g, \  (q_r, K_r, \Delta_r)\in\Acal_r,\\
& \qquad \hspace*{-12pt}\begin{rcases}
&\delta^\omega_g \in [0, \Delta_g], \\
&\displaystyle\sum_{g\in\Gfk} \delta^\omega_g = \displaystyle\sum_{r \in \Rfk} \Delta_r \ind{p^{\omega,*}\geq K_r}
\end{rcases} \quad \prob-\text{a.s.},\\
& \qquad\text{for each } \ g \in \Gfk, \ r \in \Rfk,
\end{alignedat}
\end{align}
over $\(q_g, K_g, \Delta_g \) \in \Rset^3_+$, $\Fcal$-measurable maps $\delta^\omega_g: \Omega \rightarrow [0, \Delta_g]$ for each $g \in \Gfk$, and $\(q_r,  K_r, \Delta_r \) \in \Rset^3_+$ for each $r \in \Rfk$. The function $f$ is a real-valued $\Fcal$-measurable map over all the optimization variables.
\item Buyer $r$ pays $q_r^* \Delta^*_r$ to $M$.
\item $M$ pays $q_g^* \Delta^*_g$ to seller $g$.
\item $M$ is left with a day-ahead merchandising surplus of
$$\sum_{r \in \Rfk} q^*_r \Delta^*_r - \sum_{g \in \Rfk} q^*_g \Delta^*_g.$$
\end{itemize}

\subsubsection{Real-time stage:}
\begin{itemize}[leftmargin=*]
\item Scenario $\omega$ is realized, and the real-time price of electricity $p^{\omega,*}$ is computed by $M$. 

\item $M$ pays $\left(p^{\omega,*} - K^*_r\right)^+ \Delta^*_r$ to buyer $r$. 

\item Seller $g$ pays $\left(p^{\omega, *} - K^*_g\right)^+\delta_g^{\omega, *}$ to $M$.

\item $M$ is left with a real-time merchandising surplus of
$$ -\sum_{r \in \Rfk} \left(p^{\omega,*} - K^*_r\right)^+ \Delta^*_r + \sum_{g \in \Gfk} \left(p^{\omega, *} - K^*_g\right)^+ \delta_g^{\omega, *}.$$

\end{itemize}

\subsubsection{Explaining the market clearing procedure:} The volume of call options bought by the renewable power producers are deemed to equal the volume sold by the dispatchable generators at the forward stage. The option prices and the volumes are decided in a way that each triple $\(q_i, K_i, \Delta_i \)$ lies in the set of acceptable trades $\Acal_i$ for each participant $i \in \Gfk \cup \Rfk$. In real-time, if the price $p^{\omega, *}$ exceeds the strike price $K_r$ for a renewable power producer, it will exercise its right to receive the cash amount  $\( p^{\omega, *} - K_r \)^+ \Delta_r$. The total volume of encashed call options equals the sum of $\Delta_r$ over the set of buyers for whom the spot price exceeds their strike price. $M$ allocates this volume among the option sellers. That is, it decides $\delta^\omega_g \in [0, \Delta_g]$ for each seller $g$ in a way that the volume of call options encashed equals the total volume allocated to the sellers.

\subsubsection{How market participant $i$ decides $\Acal_i$:}
Consider a seller $g \in \Gfk$ who expects a revenue $\pi_g^\omega$ in scenario $\omega$. By participating in the energy market and the option trade with the triple $\(q_g, K_g, \Delta_g \)$, she receives a payoff of 
$$ \pi_g^\omega + q_g \Delta_g - \( p^{\omega,*} - K_g \)^+ \delta_g^\omega,$$
if $M$ allocates $\delta_g^\omega \in [0, \Delta_g]$ in scenario $\omega$. Having no control over $\delta_g^\omega$, assume that seller $g$ regards $M$ as adversary. That is, she conjectures that $M$ chooses $\delta_g^\omega$ in a way that minimizes $g$'s payoff in scenario $\omega$:
$$\Pi^\omega_g \( q_g, K_g, \Delta_g \) := \pi_g^\omega + q_g \Delta_g - \( p^{\omega,*} - K_g \)^+ \Delta_g.$$
Then, a risk-neutral market participant will accept the trade defined by $\(q_g, K_g, \Delta_g \)$, if
$$ \E \[\Pi^\omega_g \( q_g, K_g, \Delta_g \) \] \geq \E \[ \pi^\omega_g \].$$
That is, $\Acal_g$ is comprised of the trades that satisfy the above inequality. Similarly, $\Acal_r$ for a risk-neutral buyer $r \in \Rfk$ is comprised of the trades defined by $\(q_r, K_r, \Delta_r \)$ that satisfy $ \E \[\Pi^\omega_r \( q_r, K_r, \Delta_r \) \] \geq \E \[ \pi^\omega_r \],$
where 
$$ \Pi^\omega_r \( q_r, K_r, \Delta_r \) := \pi^\omega_r - q_r \Delta_r + \( p^{\omega,*} - K_r \)^+ \Delta_r.$$

If a market participant is risk-averse, one can replace the expectation with an appropriate risk-functional; see Section 5 for more details.
%
%
%

\subsubsection{The objective function $f^\omega$:} The function $f^\omega$ appearing in \eqref{eq:DA.opt} depends on the nature of the market maker. As an example, consider $M$ to be a risk-neutral profit-maximizer. Then, $f^\omega$ is the merchandising surplus for $M$ in scenario $\omega$, given by  
\begin{align*}
{\sf MS}^\omega &:= \sum_{r \in \Rfk} q_r \Delta_r - \sum_{g \in \Gfk} q_g \Delta_g \\
& \qquad -\sum_{r \in \Rfk} \left(p^{\omega,*} - K_r\right)^+ \Delta_r + \sum_{g \in \Gfk} \left(p^{\omega, *} - K_g\right)^+ \delta_g^{\omega}.
\end{align*}



\subsection{Example power system with centralized option market}\label{sec:central.example}

Recall that we presented a bilateral trade of a call option between the wind power producer $W$ and a peaker power plant $P$ in the example in Section \ref{sec:motivate}. We now reconsider the same example and demonstrate how the centralized clearing mechanism generalizes the bilateral trade.

Let $W$ be the only option buyer and $P$ be the only option seller. $M$ decides  $\Acal_0$, the set of acceptable option prices, strike prices, and trade volumes. For illustrative purposes and to avoid trivial trades in our discussion, consider
$$\Acal_0 := (0, \frac{1}{\rho}]\ \times (0, \frac{1}{\rho}] \times (0, \sqrt{3}\sigma ] \subset \Rset^3_{++}.$$

We assume that $M$ caps the trade volume to equal the maximum energy shortfall in available wind from its forward contract $X_W^*$. Assuming the players to be risk-neutral, the set of acceptable trades for $P$ and $W$ are given by
\begin{align}
&\Acal_P = \{( q_P, K_P, \Delta_P)\in\Acal_0:  \nonumber \\
		& \qquad \qquad K_P+2q_P\geq 1/\rho,\Delta_P \in (0, \sqrt{3} \sigma] \}, \label{eq:AP}\\
&\Acal_W = \{( q_W, K_W, \Delta_W)\in\Acal_0: \nonumber\\
		& \qquad \qquad K_W+2q_W\leq1/\rho,\Delta_W  \in (0, \sqrt{3} \sigma] \}. \label{eq:Aw}
\end{align}

\begin{figure*}
\begin{center}
\includegraphics[width=0.6\textwidth,height=1.5in]{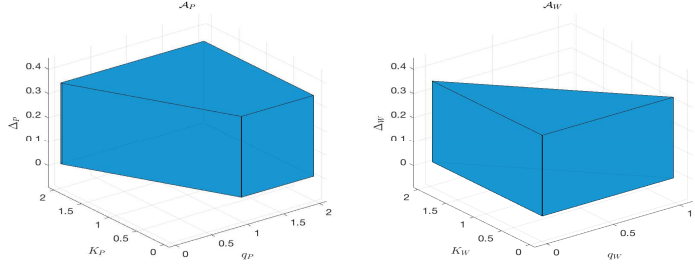}    
\caption{\small{Sets of acceptable trades $\Acal_P$ and $\Acal_W$ with $\sigma = 0.2$ and $\rho = 0.5$ are respectively drawn on the left and right. These sets are for the case when $P$ and $W$ are risk-neutral.}}
\label{fig:neutral}
\end{center}
\end{figure*}

Assuming $M$ to be risk-neutral, she solves (at day-ahead)
\begin{align}
\label{eq:example.opt}
\begin{alignedat}{5}
&{\text{maximize}} \\  
& \ \ \ \E[q_W\Delta_W-q_P\Delta_P\\
& \qquad \qquad-(p^{\omega,*} - K_W)^+ \ \Delta_W+(p^{\omega,*} - K_P)^+ \ \delta^\omega_P], \\
& \text{subject to} \\
&  \ \ \ \Delta_P = \Delta_W, (q_P, K_P, \Delta_P)\in\Acal_P, (q_W, K_W, \Delta_W)\in\Acal_W,\\
&  \ \ \ \hspace*{-12pt}\begin{rcases}
&\delta^\omega_P \in [0, \Delta_P], \\
&\displaystyle\delta^\omega_P = \Delta_W \ind{p^{\omega,*}\geq K_W}
\end{rcases} \quad \prob-\text{a.s.},\\
\end{alignedat}
\end{align}
%

We characterize the optimal solutions of \eqref{eq:example.opt} in the following result, whose proof can be found in Appendix A. 
%
\begin{prop}
\label{prop:central}
The triple $\( q_i^*, K_i^*, \Delta_i^* \) \in \Rset^3_+$ given by
$$ 2q^*_i + K^*_i = \frac{1}{\rho}\quad \text{ and } \quad \Delta^*_i \in (0,\sqrt{3}\sigma],$$
for each $i \in \lb W, P \rb$, constitutes the optimal solution of \eqref{eq:example.opt}. Moreover, at $\Delta^*_i=\sqrt{3}\sigma$, the variances of the total payments satisfy
\begin{align*}
{\sf var}[\Pi^\omega_i\( q_i^*, K_i^*, \Delta_i^* \)] - {\sf var}[\pi^\omega_i] 
= -3K^*_i\sigma^2/2< 0.
\end{align*}
\end{prop}
We note that volatility in payments for $W$ and $P$ always decrease, for any other optimal $\Delta^*_W$ and $\Delta^*_P$. However, we restrict our attention to $\Delta^*_W=\Delta^*_P=\sqrt{3}\sigma$ for illustrative purposes. Proposition 2 reveals that the centralized option market allows both $W$ and $P$ to reduce their payment volatilities. This reduction closely resembles the one in Proposition 1. The exact value of reduction, however,  depends on the strike prices that $M$ chooses. The merchandising surplus of $M$ in scenario $\omega$ is given by
\begin{align*}{ \sf MS}^{\omega, *}
= \begin{cases} 
\(q^*_P - q^*_W\)\sqrt{3}\sigma, & \text{if}\ \ \omega \leq \mu,\\
\(q^*_W-q^*_P \)\sqrt{3}\sigma, & \text{otherwise}.
\end{cases}
\end{align*}
The above expression implies that the expected merchandising surplus of the market maker remains zero, irrespective of the choice of $q_W^*$ and $q_P^*$. A risk-neutral market maker is agnostic to that choice, which however affects the volatilities of the market participants differently.
The choice $q^*_W = q^*_P$ yields the outcome for the bilateral trade, in which case the market maker only facilitates passing the payments from one player to another in each scenario $\omega$. We note that such a choice will be appropriate, if $M$ is motivated by a social goal, e.g., when $M$ is the $SO$.\footnote{By setting ${MS}^\omega$ to zero for all scenarios, the optimization problem (\ref{eq:example.opt}) can be reformulated into a system of nonlinear equations, which can be solved by the widely used Newton-Raphson's method, and it can be verified that it converges to  $q^*_W=q^*_P$.} Our discussion serves to illustrate that the bilateral trade is a special case of the centralized clearing mechanism.

Consider the case where another generator joins the market, a peaker plant $P'$ that has a linear cost of production with a marginal cost higher than that of $P$. Again, suppose $P'$ has an infinite capacity of production. It can be verified that $P'$ will never be dispatched, and hence, does not get paid from the electricity market. $P'$ can still participate in the options trade, wherein she fulfills the role of a seller. $M$ then allocates the option demand from $W$ between $P$ and $P'$. One can show that the profits and losses from the options trade in various scenarios will then be shared between $P$ and $P'$. While the option trade will still reduce the volatility of the payments for $P$, the absolute reduction in volatility will be lower. The payments to $W$, however, remain unaffected by the presence of $P'$.

\section{Generalizations}
\label{sec:risk}

\begin{figure*}
\begin{center}
\includegraphics[width=0.6\textwidth,height=1.5in]{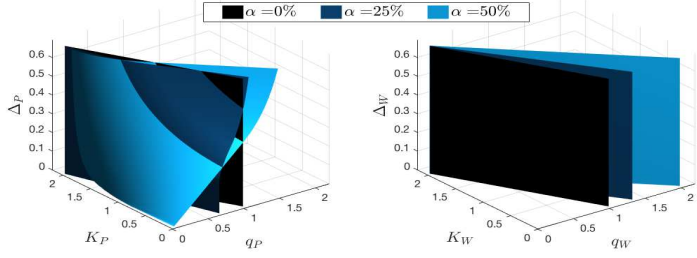}    
\caption{\small{The boundaries of the sets of acceptable trades $\Acal_P$ and $\Acal_W$ are portrayed respectively on the left and the right, for the example in Section \ref{sec:central.example} with $\mu = 1$, $\sigma = 0.4$, $\rho = 0.5$. We let both $P$ and $W$ measure risk-aversion using ${\sf CVaR}_\alpha$. The sets are computed by leveraging the technique outlined in \cite[equation (6)]{Hong}.}}
\label{fig:Aw}
\end{center}
\end{figure*}


\vspace{-0.1in}
Recall that we described the set of acceptable trades $\Acal_i$ for each $i \in \Gfk \cup \Rfk$ by assuming each market participant to be risk-neutral in Section \ref{sec:centralized}. A trade, identified by the triple $\( q_i, K_i, \Delta_i \)$, is deemed acceptable, if the expected payment from participation in the option trade together with the electricity market improves upon the expected payment from the electricity market alone. Market participants in electricity markets are known to be risk-averse in practice; e.g., see \cite{Bjorgan2,der}. Here, we model such risk-aversion and demonstrate its effect on the sets of acceptable trades.

Assume that a risk-averse market participant $i \in \Gfk \cup \Rfk$ finds a trade triple $\( q_i, K_i, \Delta_i \)$ acceptable, if 
\begin{align}
{\sf CVaR}_{\alpha_i} \[- \Pi^\omega_i \(q_i, K_i, \Delta_i \)\] \leq {\sf CVaR}_{\alpha_i} \[ -\pi^\omega_i \],
\label{eq:cvar.i}
\end{align}
where  
$$ {\sf CVaR}_\alpha \[z^\omega\] := \min_{t\in\Rset} \lb t + \frac{1}{1-\alpha} \E\[ \(z^\omega - t \)^+\]\rb,$$
$z$ is a real-valued $\Fcal$-measurable map, and $\alpha \in [0,1)$ is a parameter that quantifies the extent of risk-aversion. $\sf CVaR$ is the popular \emph{conditional value at risk} functional that describes a coherent risk measure; see \cite{risk} for a survey. If $z^\omega$ is the monetary loss in scenario $\omega$, ${\sf CVaR}_\alpha \[z^\omega \]$ equals the expected loss over the $\alpha \%$ scenarios that result in the highest losses. 

Consider the set of acceptable trades described by \eqref{eq:cvar.i} for market participant $i$, when her parameter for risk-aversion is $\alpha_i$. To make its dependence on $\alpha_i$ explicit, call this set $\Acal_i (\alpha_i)$. The risk-neutral case is modeled as $\alpha_i = 0$. The closer $\alpha_i$ is to $1$, the more risk-averse $i$ is. Notice that the set $\Acal_i(0)$ depends on the payments from the option trade alone. When $\alpha_i > 0$, the set $\Acal_i(\alpha_i)$ also depends on the payments from the electricity market.

Recall that we characterized $\Acal_W(0)$ in \eqref{eq:Aw} when $W$ is risk-neutral. To study the effect of risk-aversion, we plot the boundary of $\Acal_W(\alpha_W)$ for various values of $\alpha_W$. In these plots, we let the cap on $\Delta$ be $2\sqrt{3}/5$ to better illustrate the sets in Figure \ref{fig:Aw} (acceptable trades are on the left of the boundary). The larger the $\alpha_W$ is, the larger is the set of acceptable trades. Intuitively, a more risk-averse $W$ will accept higher option/strike prices to mitigate the real-time financial risks. In a sense, ${\sf CVaR}_{\alpha_W}$ puts more weight on monetary losses as $\alpha_W$ grows. For this example, we have 
\begin{align*} 
& \(q_W, K_W, 2\sqrt{3}/5\sigma \) \in \Acal_W\( \alpha_W \) \\
& \quad \implies \(q_W, K_W,  \Delta_W \)\in \Acal_W\( \alpha_W \) \ \forall \ \Delta_W \in (0,2\sqrt{3}/5\sigma].
\end{align*}
The boundary of $\Acal_P$ is portrayed in Figure \ref{fig:Aw}, and the acceptable set is on the right of the boundary. The changes in the sets $\Acal_P(\alpha)$ and $\Acal_W(\alpha)$ with $\alpha$ are qualitatively different. The difference stems from the fact that $W$ is exposed to negative payments while $P$ is not. In a sense, more risk aversion makes $W$ more prone to purchase the options, while $P$ becomes less likely to sell more options.


%
%
 
\subsubsection{Effect of incorrect price/payment conjectures:}
\vspace{-0.1in}

If a market participant, say $W$, has an incorrect price conjecture, denoted by $\hat{p}^{\omega,*}$, then her set of acceptable trades $\hat{\Acal}_W$ can be different from $\Acal_W$. While we assumed in our examples that participants have correct price conjectures, our model can be generalized to the case where a seller or a buyer bids an incorrect set. Furthermore, even if $W$ has the correct conjecture, she might have an incentive to misrepresent her acceptability set. Such considerations are beyond the scope of this paper and can be investigated in future work. 
%



\section{Conclusions}

\label{sec:conc}
This paper has explored the design of a centralized cash-settled options market that can be implemented in parallel with an existing electricity market. Through stylized examples, we have demonstrated that our design leads to reduction in the volatility in payments. Additionally, we have shown that a renewable supplier is less exposed to negative payments under uncertainty, compared to a conventional approach. The paper serves as a stepping stone for a more comprehensive market with deep integration of renewable supply. An attractive feature of the proposed design is that it can be applied to any existing electricity market, making its implementation less complex. Generalizations include: effect of risk-aversion, effect of price conjectures, computational aspects, and multi-period considerations.

\appendix

\section{Proofs of Propositions}

\vspace{-0.1in}

\subsection{Proof of Proposition 1}
\vspace{-0.1in}

With a slight abuse of notation, we sometimes use $\Delta$ for $\Delta(q,K)$. Suppose $P$ chooses the pair $\( q, K \) \in \Rset^2_+$. Then, the payoff of $W$ from the option trade alone is given by
\begin{align*}
V^\omega_W(q, K, \Delta)
:= \Pi^\omega_W (q, K, \Delta) - \pi^\omega_W. 
\end{align*}
Utilizing the relation in \eqref{eq:PiW}, we have
\begin{align}
& \E\[V^\omega_W(q, K, \Delta)\] \notag\\
& \qquad= \begin{cases}
-q\Delta, & \text{if } K > 1/\rho,\\
- \frac{\Delta}{2} \(  2q +  K - 1/{\rho} \), & \text{otherwise}.
\end{cases}
\label{eq:evw}
\end{align}
We characterize the best response of $W$ to the choice of $\(q,K\)$ by $P$. Call it the set-valued map $\Delta^*(q, K)$. If $K > 1/\rho$, then $\Delta^*(q, K) = \{ 0 \}$. Otherwise,
    \begin{itemize}
    \item If $2q + K < 1/\rho$, then $\Delta^*(q, K) = \{ \sqrt{3}\sigma \}$.
    \item If $2q + K = 1/\rho$, then $W$ is agnostic to the choice of $\Delta$, and hence, $\Delta^*(q, K) = [0, \sqrt{3}\sigma]$.
    \item If $2q + K >1/\rho$, then  $\Delta^*(q, K) = \{ 0 \}$.
    \end{itemize}

Similarly, define $V^\omega_P(q, K, \Delta):= \Pi^\omega_P (q, K, \Delta) - \pi^\omega_P,$
as the payoff of $P$ from the option trade. Then, \eqref{eq:PiW} and \eqref{eq:PiP} imply that
\begin{align}
\E \[ V^\omega_P(q, K, \Delta) \] = - \E \[ V^\omega_W(q, K, \Delta) \].
\label{eq:evp}
\end{align}
Given the best response of $W$, we have the following cases.  
  \begin{itemize}
    \item If $2q + K < 1/\rho$, then $\E \[ V^\omega_P(q, K, \Delta) \]<0$, and hence no trades will occur and $P$ must choose $(q,K)$ such that $2q + K \geq 1/\rho$.
    \item If $2q + K = 1/\rho$, then $\E \[ V^\omega_P(q, K, \Delta) \]=0$, and with $\Delta^*(q, K) = [0, \sqrt{3}\sigma]$, this is the only case at which nontrivial trades occur.
    \item If $2q + K >1/\rho$, and with $\Delta^*(q, K) = \{ 0 \}$, trivial trades occur.
    \end{itemize}

Thus, the Stackelberg equilibria of $\Gcal$ are given by $\Ncal_1 \cup \Ncal_2$, where
\begin{align*}
\Ncal_1 &:= \left\{ (q, K, \Delta) \ \vert \ (q,K) \in \Rset_+^2, \ \Delta:\Rset^2_+ \to [0, \sqrt{3}\sigma], \right.\\
&\qquad\quad \left. 2q + K > 1/\rho , \ \Delta(q', K') = 0 \ \forall \ (q', K') \in\Rset^2_+  \right\},\\ 
\Ncal_2 &:= \left\{ (q, K, \Delta) \ \vert \ (q,K) \in \Rset_+^2, \ \Delta:\Rset^2_+ \to [0, \sqrt{3}\sigma], \right.\\ &\qquad\quad \left. 2q + K = 1/\rho \right\}.
\end{align*}
 For any Stackelberg equilibrium $(q^*, K^*, \Delta^*(q^*,K^*)=\sqrt{3}\sigma)\in\Ncal_2$, we have $V^{\omega,*}_W= q^*\sqrt{3}\sigma$, if $\omega<\mu$. Otherwise, $V^{\omega,*}_W=- q^*\sqrt{3}\sigma$. The differences in variances are given by \begin{eqnarray*}
&{\sf var}\left[\Pi^\omega_W(q^*, K^*, \Delta^*(q^*,K^*))\right]- {\sf var}\left[\pi^\omega_W\right]\\
&\qquad \qquad =2{\sf cov}(\pi^\omega_W,V^{\omega,*}_W)+{\sf var}\left[V^{\omega,*}_W\right], \end{eqnarray*}
where $\pi_W^\omega = \mu -( \mu - \omega)^+/\rho$. Using the fact that $1/\rho =2q^*+K^* $, it can be verified upon substitutions that  $${\sf var}\left[\Pi^\omega_W(q^*, K^*, \Delta^*(q^*, K^*))\right]- {\sf var}\left[\pi^\omega_W\right]=-3K^*\sigma^2/2 <0.$$
For $P$, we have $\pi_P^\omega=\mu-\pi_W^\omega \ \ \text{and} \ \ V^{\omega,*}_P=-V^{\omega,*}_W,$ and the statement of the Proposition can be easily verified.

\subsection{Proof of Proposition 2}
\vspace{-0.15in}
The constraints in (\ref{eq:example.opt}) imply that $\Delta_P=\Delta_W$. By (\ref{spot}), it can be checked that 
\begin{align*}
&\E[{ \sf MS}^\omega]=(2q_W+K_W-1/\rho)\frac{\Delta_W}{2}\\
&\qquad \qquad +(-2q_P-K_P+1/\rho)\frac{\Delta_P}{2},
\end{align*}
Given $\Acal_P$ described in (\ref{eq:AP}), $M$ must choose $2q_P+K_P=1/\rho$, and it will be agnostic to $\Delta_P$. Similarly for $W$, $M$ must choose $2q_W+K_W=1/\rho$, and hence the maximum is achieved at 
$$ 2q^*_i + K^*_i =1/\rho\quad \text{ and } \quad \Delta_P=\Delta_W \in (0,\sqrt{3}\sigma],$$
where $i \in \lb W, P \rb$. The proof of the variances follows the exact steps in the proof of Proposition $\ref{prop:bilat}$.

\vspace{-0.1in}

\bibliography{ifacconf}

\end{document}